# Stability Analysis of Picard Iteration for Coupled Neutronics/Thermal-Hydraulics Simulations

Dean Wang

*Nuclear Engineering Program, The Ohio State University, Columbus, OH 43210*
*wang.12239@osu.edu*

## INTRODUCTION

Reactor core analysis often needs to solve a multiphysics nonlinear coupled system, including neutron transport, thermal-hydraulics, and other important physics phenomena. One straightforward method for solving such a coupled system is Picard fixed-point iteration [1], which alternates between solving individual physics problems separately. However, many numerical studies show that Picard iteration can be unstable, and a user-defined relaxation is usually required to achieve convergence [2-4].

In this paper, we present a formal Fourier analysis (FA) of Picard iteration for the coupled neutronics/thermal hydraulics (N/TH) problem and derive theoretical predictions for the spectral radius of Picard iteration for such coupled calculations as a function of the temperature difference between the fuel and coolant, temperature coefficients of cross sections (i.e., Doppler feedback), scattering ratio, and core height. An optimal underrelaxation factor is also derived based on the Fourier analysis.

## FORMULATION AND ALGORITHM

We consider the following simple one-group, planar-geometry k-eigenvalue problem on the domain $0 \leq x \leq L$ with reflective boundary conditions:

$$\mu \frac{\partial \psi(x,\mu)}{\partial x} + \Sigma_t(T)\psi(x,\mu) = \frac{1}{2}\Sigma_s(T)\phi(x) + \frac{1}{2k_{eff}}\nu\Sigma_f(T)\phi(x), \quad (1)$$

and the simplified heat transfer equation for a single typical pressurized water reactor (PWR) fuel pin:

$$T = T_m + A\Sigma_f(T)\phi(x), \quad (2)$$

with

$$A = \pi r_{fo}^2 \kappa R_t, \quad (3a)$$

and

$$R_t = \left[\frac{1}{8\pi k_f} + \frac{1}{2\pi r_g h_g} + \frac{1}{2\pi k_c}\ln\left(\frac{r_{co}}{r_{ci}}\right) + \frac{1}{2\pi r_{co} h}\right], \quad (3b)$$

where

$\psi$ = neutron angular flux
$\phi = \int_{-1}^{1} \psi(x,\mu)d\mu$, neutron scalar flux
$\Sigma_t$ = macroscopic total cross section
$\Sigma_s$ = macroscopic scattering cross section
$\Sigma_f$ = macroscopic fission cross section
$\nu$ = average neutron yield per fission
$k_{eff}$ = effective multiplication factor
$\kappa$ = average energy released per fission
$T$ = volume averaged fuel temperature
$T_m$ = bulk coolant temperature
$r_{fo}$ = fuel radius
$r_{ci}$ = cladding inner radius
$r_{co}$ = cladding outer radius
$r_g = \frac{r_{ci}+r_{fo}}{2}$, mean radius in the gap
$k_f$ = fuel thermal conductivity
$k_c$ = cladding thermal conductivity
$h_g$ = effective gap conductance
$h$ = coolant convection heat transfer coefficient

Note that the linear heat generation rate (or linear power) of the fuel rod, $q'$, can be calculated by

$$q'(x) = \pi r_{fo}^2 \kappa \Sigma_f(T)\phi(x). \quad (4)$$

Picard iteration is used to solve the above coupled N/TH system as follows. The transport equation is solved first, then the fuel temperature is calculated using the newly obtained thermal power (neutron flux). An underrelaxation factor is introduced in the temperature update. Note that here the transport iteration is fully converged during each TH update. This type of coupling is also referred to as the "loose" coupling method, in comparison to the "tight" coupling where the TH update is performed during transport sweeping.

$$\mu \frac{\partial \psi^{(k+1)}(x,\mu)}{\partial x} + \Sigma_t(T^{(k)})\psi^{(k+1)} = \frac{1}{2}\Sigma_s(T^{(k)})\phi^{(k+1)}(x) + \frac{1}{2k_{eff}}\nu\Sigma_f(T^{(k)})\phi^{(k+1)}, \quad (5)$$

$$T^* = T_m + A\Sigma_f(T^{(k)})\phi^{(k+1)}(x), \quad (6a)$$

$$T^{(k+1)} = \omega T^* + (1-\omega)T^{(k)}, \quad (6b)$$

where $\omega$ is the underrelaxation factor and the superscript $k$ denotes the iteration number.

## LINEARIZATION

To perform Fourier analysis of the coupled N/TH problem, we need to first linearize the system of equations. We define the following linearized variables:

$$\psi(x,\mu) = \psi_0(x,\mu) + \varepsilon\psi_1(x,\mu), \quad (7a)$$

$$\phi(x) = \phi_0(x) + \varepsilon\phi_1(x), \quad (7b)$$

$$k_{eff} = k_{eff,0}, \tag{7c}$$

$$T(x) = T_0 + \varepsilon T_1(x), \tag{7d}$$

$$\Sigma_i(T) = \Sigma_{i0} + \Sigma_{i1}(T - T_0)$$
$$= \Sigma_{i0} + \varepsilon \Sigma_{i1} T_1(x), \quad i = t, s, f, a \tag{7e}$$

Note that $k_{eff} = k_{eff,o}$ due to the flux normalization. The cross sections are assumed to be linearly dependent on the fuel temperature. However, other feedback mechanisms such as thermal expansion [5] and moderator temperature feedbacks can be treated as well.

Substituting Eqs. (7a) - (7e) into (5), after some algebra we obtain by neglecting the $O(\varepsilon^2)$ terms

$$\mu \frac{\partial \psi_1^{(k+1)}(x,\mu)}{\partial x} + \Sigma_{t0}\psi_1^{(k+1)}(x,\mu) + \Sigma_{t1}T_1^{(k)}(x)\psi_0(x,\mu)$$
$$= \frac{1}{2}\Sigma_{s0}\phi_1^{(k+1)}(x) + \frac{1}{2}\Sigma_{s1}T_1^{(k)}(x)\phi_0$$
$$+ \frac{1}{2k_{eff,0}}\nu\Sigma_{f0}\phi_1^{(k+1)}(x) + \frac{1}{2k_{eff,0}}\nu\Sigma_{f1}T_1^{(k)}(x)\phi_0. \tag{8}$$

For reflective BC, $\psi_0 = \frac{\phi_0}{2}$, and $\Sigma_{a0} = \frac{1}{k_{eff,0}}\nu\Sigma_{f0}$, then we rewrite Eq. (8) as

$$\mu \frac{\partial \psi_1^{(k+1)}(x,\mu)}{\partial x} + \Sigma_{t0}\psi_1^{(k+1)}(x,\mu)$$
$$= \frac{1}{2}\Sigma_{t0}\phi_1^{(k+1)}(x) - \frac{1}{2}\Sigma_{t0}\gamma T_1^{(k)}(x), \tag{9}$$

where

$$\gamma = (1-c_0)\left(\frac{\Sigma_{a1}}{\Sigma_{a0}} - \frac{\Sigma_{f1}}{\Sigma_{f0}}\right)\phi_0, \tag{10a}$$

with

$$\Sigma_{a1} = \Sigma_{t1} - \Sigma_{s1}, \tag{10b}$$

$$c_0 = \frac{\Sigma_{s0}}{\Sigma_{t0}}. \tag{10c}$$

Substituting Eqs. (7b), (7d), and (7e) into (6a) and (6b) respectively, we obtain

$$T_1^*(x) = A\Sigma_{f0}\phi_1^{(k+1)}(x) + A\Sigma_{f1}\phi_0 T_1^{(k)}(x), \tag{11}$$

$$T_1^{(k+1)}(x) = \omega T_1^*(x) + (1-\omega)T_1^{(k)}(x). \tag{12}$$

Then we substitute Eq. (11) into (12) to give

$$T_1^{(k+1)}(x)$$
$$= \omega A\Sigma_{f0}\phi_1^{(k+1)}(x) + (1 - \omega + \omega A\Sigma_{f1}\phi_0)T_1^{(k)}(x). \tag{13}$$

For brevity we drop the subscript "1" in the flux and temperature variables without confusion

$$\mu \frac{\partial \psi^{(k+1)}(x,\mu)}{\partial x} + \Sigma_{t0}\psi^{(k+1)}(x,\mu)$$
$$= \frac{1}{2}\Sigma_{t0}\phi^{(k+1)}(x) - \frac{1}{2}\Sigma_{t0}\gamma T^{(k)}(x), \tag{14}$$

$$T^{(k+1)}(x) =$$

$$\omega A\Sigma_{f0}\phi^{(k+1)}(x) + (1 - \omega + \omega A\Sigma_{f1}\phi_0)T^{(k)}(x). \tag{15}$$

**FOURIER ANALYSIS**

We introduce the inverse Fourier transforms:

$$\phi^{(k)}(x) = \int_{-\infty}^{+\infty} a^{(k)}(\xi)e^{i\Sigma_{t0}\xi x}d\xi, \tag{16a}$$

$$\psi^{(k)}(x,\mu) = \int_{-\infty}^{+\infty} b^{(k)}(\xi,\mu)e^{i\Sigma_{t0}\xi x}d\xi, \tag{16b}$$

$$T^{(k)}(x) = \int_{-\infty}^{+\infty} c^{(k)}(\xi)e^{i\Sigma_{t0}\xi x}d\xi. \tag{16c}$$

The same Fourier ansatz is used for the temperature as for the neutron flux because the fuel temperature is roughly proportional to the neutron flux as shown in Eq. (6a). The solutions are required to satisfy the boundary conditions. The discrete Fourier error mode $\xi$ for the reflective boundary conditions are given below in Eq. (17). If the periodic boundary conditions are used, then they are simply multiplied by a factor of 2.

$$\xi = \frac{\pi}{\Sigma_{t0}L}j, \quad j = \pm 1, \pm 2, \ldots \tag{17}$$

where $L$ is the reactor core height (or the fuel rod length), i.e., the slab thickness in our model problem.

By substituting Eqs. (16a) - (16c) into (14) and noting that each of the Fourier modes is independent, we obtain

$$\Sigma_{t0}(i\xi\mu + 1)b^{(k+1)}(\xi,\mu)$$
$$= \frac{1}{2}\Sigma_{t0}a^{(k+1)}(\xi) - \frac{1}{2}\Sigma_{t0}\gamma c^{(k)}(\xi). \tag{18}$$

We rewrite Eq. (18) as

$$b^{(k+1)}(\xi,\mu) = \frac{1}{2}\frac{1}{(i\xi\mu+1)}a^{(k+1)}(\xi) - \frac{1}{2}\frac{\gamma}{(i\xi\mu+1)}c^{(k)}(\xi). \tag{19}$$

By noting that $\int_{-1}^{1}\frac{1}{2}\frac{1}{(i\xi\mu+1)}d\mu = \tan^{-1}(\xi)/\xi$ and $a^{(k+1)}(\xi) = \int_{-1}^{1}b^{(k+1)}(\xi,\mu)d\mu$, we integrate the above equation with respect to $\mu$ to obtain

$$a^{(k+1)}(\xi) = -\frac{\gamma\rho_{PI}(\xi)}{1-\rho_{PI}(\xi)}c^{(k)}(\xi), \tag{20}$$

where

$$\rho_{PI}(\xi) = \frac{\tan^{-1}(\xi)}{\xi}. \tag{21}$$

Note that $\rho_{PI}$ is the spectral radius function for the standard power iteration (PI) algorithm.

Substituting Eqs. (16a) and (16c) into (15), we obtain

$$c^{(k+1)}(\xi) = \omega A\Sigma_{f0}a^{(k+1)}(\xi)$$
$$+ (1 - \omega + \omega A\Sigma_{f1}\phi_0)c^{(k)}(\xi). \tag{22}$$

Substituting Eq. (20) into (22), we have

$$c^{(k+1)}(\xi) =$$

$$\left[1 - \omega + \omega A\Sigma_{f1}\phi_0 - \omega A\Sigma_{f0}\frac{\gamma\rho_{PI}(\xi)}{1-\rho_{PI}(\xi)}\right]c^{(k)}(\xi). \tag{23}$$

Thus, we obtain the spectral radius function of Picard iteration for the coupled N/TH problem as

$$\varrho(\xi) = 1 - \omega + \omega A\Sigma_{f1}\phi_0 - \omega A\Sigma_{f0}\frac{\gamma\rho_{\text{PI}}(\xi)}{1-\rho_{\text{PI}}(\xi)}. \quad (24)$$

Substituting Eqs. (3a) and (10a) into (24), we obtain

$$\varrho(\xi) = 1 - \omega \left[\begin{array}{c} 1 - \pi r_{fo}^2 \kappa R_t \Sigma_{f1}\phi_0 + \\ \pi r_{fo}^2 \kappa R_t \Sigma_{f0}\phi_0 \frac{(1-c_0)\left(\frac{\Sigma_{a1}}{\Sigma_{a0}} - \frac{\Sigma_{f1}}{\Sigma_{f0}}\right)}{1-\rho_{\text{PI}}(\xi)}\rho_{\text{PI}}(\xi) \end{array}\right]. \quad (25)$$

Substituting Eq. (4) into (25), we have

$$\varrho(\xi) = 1 - \omega\left[1 - q_0' R_t \left(\frac{\Sigma_{f1}}{\Sigma_{f0}} - \left(\frac{\Sigma_{a1}}{\Sigma_{a0}} - \frac{\Sigma_{f1}}{\Sigma_{f0}}\right)\frac{1-c_0}{1-\rho_{\text{PI}}(\xi)}\rho_{\text{PI}}(\xi)\right)\right]. \quad (26)$$

By noting that $q_0' R_t = T_0 - T_m$, Eq. (26) can be rewritten as

$$\varrho(\xi) = 1 - \omega\left\{1 - (T_0 - T_m)\left[\frac{\Sigma_{f1}}{\Sigma_{f0}} - \left(\frac{\Sigma_{a1}}{\Sigma_{a0}} - \frac{\Sigma_{f1}}{\Sigma_{f0}}\right)\frac{1-c_0}{1-\rho_{\text{PI}}(\xi)}\rho_{\text{PI}}(\xi)\right]\right\} \quad (27)$$

Finally, the spectral radius of Picard iteration for the coupled N/TH nonlinear system is given as

$$\rho = \max_\xi |\varrho(\xi)|. \quad (28)$$

**RESULTS**

The spectral radius of the Picard iteration method for the coupled N/TH system is a function of the temperature difference between the fuel and coolant, temperature coefficients of fission and absorption cross sections, scattering ratio, and spectral radius of the standard PI algorithm (or essentially the error mode).

The spectral radius function of the PI algorithm, $\rho_{\text{PI}}(\xi) = \tan^{-1}(\xi)/\xi$, attains the largest value at the error mode $\xi = \pi/(\Sigma_{t0}L)$, which is the most slowly converging mode. It is well known that the PI becomes increasingly slow ($\rho \to 1$) as the problem domain becomes large, though the method is unconditionally stable because its spectral radius always remains below 1. On the other hand, $\rho_{\text{PI}}(\xi)$ tends to zero as $\xi$ limits to infinity. In addition, it is interesting to point out that the term $\rho_{\text{PI}}(\xi)/(1 - \rho_{\text{PI}}(\xi))$ in Eq. (26) or (27) can be approximated by $3/\xi^2$ for $\xi$ small (e.g., the relative difference is less than 1% when $\xi < 0.2$). With such approximation, we have actually obtained the spectral radius for the diffusion solution coupled with TH.

For light water reactors, the cross-section temperature coefficients are typically very small. Table I summarizes the one-group cross section data for a typical PWR.

TABLE I. Typical PWR Data

| $\Sigma_{t0}$ (cm$^{-1}$) | $\nu\Sigma_{fo}$ (cm$^{-1}$) | $c_o$ | $\Sigma_{f1}/\Sigma_{f0}$ (K$^{-1}$) | $\Sigma_{a1}/\Sigma_{a0}$ (K$^{-1}$) |
|---|---|---|---|---|
| 0.718 | 0.0297 | 0.96 | $-1.99 \times 10^{-5}$ | $8.67 \times 10^{-6}$ |

Note that the temperature coefficient of the absorption cross section is positive, whereas that of the fission cross section is negaive. For such problems, the convergence of the unrelaxed Picard iteration is determined by the smallest error mode $\xi = \pi/(\Sigma_{t0}L)$, and the spectral radius is given as

$$\rho = (T_0 - T_m)\left[\left(\frac{\Sigma_{a1}}{\Sigma_{a0}} - \frac{\Sigma_{f1}}{\Sigma_{f0}}\right)\frac{1-c_0}{1-\rho_{\text{PI}}(\xi)}\rho_{\text{PI}}(\xi) - \frac{\Sigma_{f1}}{\Sigma_{f0}}\right]. \quad (29)$$

Fig. 1 shows that the spectral radius increases with the increasing reactor core height and eventually Picard iteration fails to converge when the reactor core height is larger than a critical value (for the given total cross section). If the temperature difference between the fuel and coolant increases (e.g., by increasing thermal resistance or linear heat generation rate $q'$), then the coupling becomes less stable as the spectral radius becomes larger.

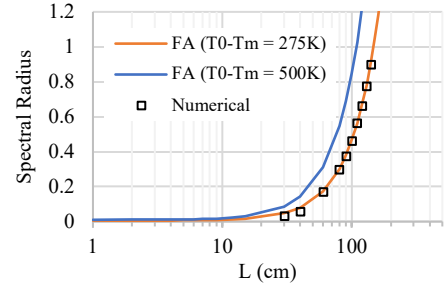

Fig. 1. Spectral radius vs. core height.

To verify the FA results, we compute numerical convergence rates based on a 1-D model problem, which is the homogeneous slab with the reflective boundary on both sides. The Gauss-Legendre S12 quadrature set is used for angular discretization and the Diamond Difference (DD) method is employed for spatial discretization. Note that the angular quadrature and the mesh size used are sufficiently fine to minimize the numerical errors. A simple heat balance model is used to calculate the fuel and coolant temperatures at each axial cell. Fig. 1 shows that the numerical results for the problem are in excellent agreement with the FA results.

For the relaxed case, i.e., $0 < \omega < 1$, underrelaxation helps to stabilize Picard iteration by improving the convergence rate as shown in Fig. 2.

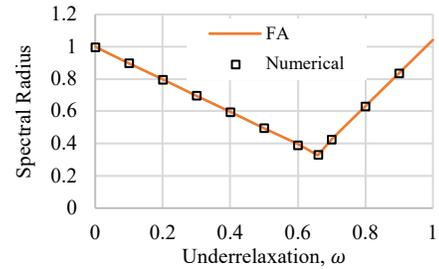

Fig. 2. Spectral radius vs. underrelaxation.

For this case, the core height $L = 150$ cm, and the typical PWR data in Table I is used. Again, the FA predictions are consistent with the numerical results.

To derive the optimal underrelaxation factor $\omega_{opt}$, it is noted that when $\omega < \omega_{opt}$, $\max_\xi |\varrho(\xi)|$ is found at $\xi \to \infty$, where $\rho_{\text{PI}}(\xi) = 0$, and

$$\rho = 1 - \omega \left[1 - (T_0 - T_m)\frac{\Sigma_{f1}}{\Sigma_{f0}}\right], \quad (30)$$

while for $\omega > \omega_{opt}$, $\max_\xi |\varrho(\xi)|$ is found at $\xi = \pi/(\Sigma_{t0}L)$, and

$$\rho = -1 + \omega \left\{1 - (T_0 - T_m)\left[\frac{\Sigma_{f1}}{\Sigma_{f0}} - \left(\frac{\Sigma_{a1}}{\Sigma_{a0}} - \frac{\Sigma_{f1}}{\Sigma_{f0}}\right)\frac{1-c_0}{1-\rho_{\text{PI}}\left(\frac{\pi}{\Sigma_{t0}L}\right)}\rho_{\text{PI}}\left(\frac{\pi}{\Sigma_{t0}L}\right)\right]\right\} \quad (31)$$

Then the optimal $\omega_{opt}$ can be obtained by equating Eqs. (30) and (31):

$$\omega_{opt} = \frac{2}{2 - (T_0 - T_m)\left[2\frac{\Sigma_{f1}}{\Sigma_{f0}} - \left(\frac{\Sigma_{a1}}{\Sigma_{a0}} - \frac{\Sigma_{f1}}{\Sigma_{f0}}\right)\frac{1-c_0}{1-\rho_{\text{PI}}\left(\frac{\pi}{\Sigma_{t0}L}\right)}\rho_{\text{PI}}\left(\frac{\pi}{\Sigma_{t0}L}\right)\right]}. \quad (32)$$

For the case shown in Fig 2, the FA predicted optimal underrelaxation factor is the same as the numerical result, $\omega_{opt} = 0.66$. Note that this case is unstable ($\rho = 1.042$) unless underrelaxation is applied. It indicates that the theoretical estimate of the optimal underrelaxation factor is quite accurate. The optimal underrelaxation depends on various parameters as indicated by Eq. (32). For example, it varies with the core height (for this case, $\Sigma_{t0} = 0.718$ cm$^{-1}$) as depicted in Fig. 3. The more underrelaxation is needed for higher cores (or longer fuel rods). It also indicates that the higher fuel/coolant temperature difference (i.e., larger linear power or thermal resistance), the more underrelaxation is necessitated for stabilizing Picard iteration.

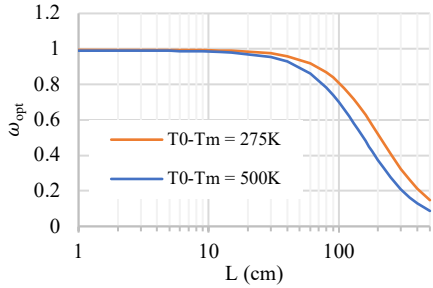

Fig. 3. Optimal underrelaxation vs. core height.

## CONCLUSIONS

We have presented a formal Fourier analysis to theoretically predict the convergence properties of Picard fixed-point iteration for coupled neutronics/thermal-hydraulics calculations. The work provides a more rigorous theoretical basis for applications of Picard iteration for such calculations. The derived closed form estimate for the spectral radius of the Picard coupling method is a function of various reactor parameters such as the fuel and coolant temperature difference (which instead depends on the rod linear power and thermal resistance), fuel temperature feedback (Doppler effect), scattering ratio, and reactor core height. It implies that Picard iteration is more stable for smaller reactors and lower rod linear power (or thermal resistance) as expected. In addition, it is worth noting that for LWRs the Doppler feedback plays a more dominant role in Picard iteration than the moderator temperature (density) feedback. This finding is consistent with numerical experiments reported in Ref. 4. We will report the analysis in the future.

A long-standing issue with Picard iteration is that it oftentimes relies on underrelaxation to stabilize the coupled calculation. However, a priori optimal underrelaxation was previously not available for a specific problem. We hope that our new theoretical result can provide a valuable estimate of underrelaxation for stabilizing coupled neutronics/thermal-hydraulics calculations.

The relaxed Picard iteration is similar to the undamped Anderson acceleration with depth $m = 1$ (AA-1), in which only one previous iterate is used [6,2,3]. However, it is expected that the Picard with optimal underrelaxation will outperform the AA-1 algorithm since the linear coefficients of AA-1 are determined by minimizing the norm of an affine combination of residual vectors and they are generally different from the optimal underrelaxation factor.

Although we have only focused on the Fourier analysis for the simple PWR model problem, the analysis presented should be applicable for other types of reactors and more realistic problems such as multigroup problems. In addition, it can be also applied to other coupling methods.